\documentclass{elsart3-1}

\usepackage{amssymb}

\usepackage[english,francais]{babel}
\usepackage{graphicx}
\usepackage{wrapfig}
\usepackage{psfrag}

\newtheorem{theorem}{Theorem}[section]
\newtheorem{lemma}[theorem]{Lemma}
\newtheorem{e-proposition}[theorem]{Proposition}

\newtheorem{e-definition}[theorem]{Definition\rm}

\newtheorem{theoreme}{Th\'eor\`eme}

\renewcommand{\theequation}{\arabic{equation}}
\def\og{\leavevmode\raise.3ex\hbox{$\scriptscriptstyle\langle\!\langle$~}}
\def\fg{\leavevmode\raise.3ex\hbox{~$\!\scriptscriptstyle\,\rangle\!\rangle$}}

\begin{document}
\begin{center}
Differential Geometry/Dynamical Systems
\end{center}
\begin{frontmatter}


\selectlanguage{english}
\title{Bifurcation and forced symmetry breaking in Hamiltonian systems}




\selectlanguage{english}
\author[authorlabel1]{F\'ethi Grabsi}
\ead{fethi.grabsi@inln.cnrs.fr}
\author[authorlabel2]{James Montaldi}
\ead{j.montaldi@umist.ac.uk}
\author[authorlabel3]{Juan-Pablo Ortega}
\ead{Juan-Pablo.Ortega@math.univ-fcomte.fr}

\address[authorlabel1]{Institut Non Lin\'eaire de Nice, UMR 129 CNRS-UNSA,
1361 route des \nobreak{Lucioles}, 06560 Valbonne, France}
\address[authorlabel2]{Department of Mathematics. UMIST. PO Box 88,
Manchester M60 1QD, United Kingdom}
\address[authorlabel3]{Centre National de la Recherche
Scientifique, D\'epartement de Math\'ematiques de Besan\c con,
Universit\'e de Franche-Comt\'e. 16, route de Gray,
25030 Besan\c con cedex, France}

\begin{abstract}
We consider the phenomenon of forced symmetry breaking in a symmetric Hamiltonian
system on a symplectic manifold. In particular we study the persistence of an
initial relative equilibrium subjected to this forced symmetry breaking. We see
that, under certain nondegeneracy conditions, an estimate can be made on the number
of bifurcating relative equilibria . {\it To cite this article: F. Grabsi et al., C.
R. Acad. Sci. Paris, Ser. I 338 (2004).}

\vskip 0.5\baselineskip

\selectlanguage{francais}
\noindent{\bf R\'esum\'e} \vskip 0.5\baselineskip \noindent {\bf Bifurcation et
brisure forc\'ee de sym\'etrie dans les syst\`emes hamiltoniens. } Nous
consid\'erons  le ph\'enom\`ene de brisure forc\'ee de sym\'etrie dans un syst\`eme
hamiltonien sym\'etrique d\'efini sur une vari\'et\'e symplectique. Plus
pr\'ecisement, nous \'etudions la persistance d'un \'equilibre relatif soumis \`a
une brisure  de sym\'etrie. Nous verrons que, sous certaines hypoth\`eses de
non-d\'eg\'en\'er\'escence, on peut donner une estimation du nombre d'\'equilibres
relatifs persistants apr\`es la brisure. {\it Pour citer cet article~: F. Grabsi et
al., C. R. Acad. Sci. Paris, Ser. I 338 (2004).}

\end{abstract}
\end{frontmatter}

\selectlanguage{francais}
\section*{Version fran\c{c}aise abr\'eg\'ee}
Nous nous int\'eressons au ph\'enom\`ene de la brisure forc\'ee de sym\'etrie dans
les syst\`emes hamiltoniens sym\'etriques. On consid\`ere  une action libre du tore
${\mathbb T}^n$ de dimension $n$  sur une vari\'et\'e symplectique $({\mathcal M},
\omega)$ qui admet une application moment (n\'ecessairement invariante) ${\mathbf
J}_{{\mathbb T}^n} : {\mathcal M} \longrightarrow (\mathfrak{t} ^n)^\ast \simeq
\mathbb{R} ^n$. Soit $H _0  $ une fonction hamiltonienne invariante par rapport \`a
l'action de ${\mathbb T}^n$ et dont le champ de vecteurs $X_{H _0}$ associ\'e
pr\'esente, par hypoth\`ese, un \'equilibre relatif (\textsc{er}). Nous rappelons
qu'un \textsc{er} $m$ de $X_{H_0}$ est d\'efini par la condition $
X_{H_0}(m)=\xi_{{\mathcal M}}(m)$, pour un \'el\'ement $\xi$ dans
l'alg\`ebre de Lie  de $\mathbb{T} ^n $ d\'enomm\'e la vitesse de
l'\textsc{er}.

Soit $ {{\mathbb T}}^r\subset {{\mathbb T}}^n$ un sous-tore de ${{\mathbb T}}^n$ et
$H_{\varepsilon}$ une famille de perturbations ${{\mathbb T}}^r$-invariantes de
l'hamiltonien $H_0$, param\'etr\'ee d'une mani\`ere lisse par $\varepsilon\in
\mathbb{R}$. En principe, les \textsc{er} du syst\`eme $({\mathcal M}, \omega, H_0,
{{\mathbb T}}^n, {\bf J}_{{{\mathbb T}}^n} )$ ne seront plus en g\'en\'eral des
\textsc{er} pour le syst\`eme $({\mathcal M}, \omega, H_{\varepsilon}, {{\mathbb
T}}^r, {\bf J}_{{{\mathbb T}}^r})$. Notre int\'er\^et est de d\'eterminer sous
quelles conditions l'\textsc{er} de d\'epart $m$ continue \`a \^etre un \textsc{er}
du syst\`eme $({\mathcal M}, \omega, H_{\varepsilon}, {{\mathbb T}}^r, {\bf
J}_{{{\mathbb T}}^r})$. Plus pr\'ecisement, nous montrons que  les \textsc{er} qui
sont non-d\'eg\'en\'er\'es dans un certain sens et dont la vitesse appartient \`a
l'alg\`ebre de Lie du sous-tore $ \mathbb{T}^r$ persistent. De plus, le th\'eor\`eme
suivant nous donne une estimation de leur nombre  :

\begin{theoreme}
Soit  $H_\varepsilon$ une famille de fonctions hamiltoniennes param\'etris\'ee d'une
mani\`ere lisse par $\varepsilon \in \mathbb{R}$ et d\'efinie sur la vari\'et\'e
symplectique $({\mathcal M},\omega)$. Supposons que $H_0$ est invariant par rapport
\`a une action libre et canonique du tore $\mathbb{T}^n$ et que le champ de vecteurs
associ\'e $X_{H_0} $ a un \'equilibre relatif $m$ avec vitesse $\xi \in
\mathfrak{t}^n $ et moment $\mu:=\mathbf{J}_{\mathbb{T}^n}(m)\in (\mathfrak{t}^n)
^\ast \simeq \mathbb{R}^n$. $\mathbf{J}_{\mathbb{T}^n}:\mathcal{M}\rightarrow
(\mathfrak{t} ^n) ^\ast $ est une application moment associ\'ee \`a l'action de
$\mathbb{T}^n$.

Soit $\mathbb{T}^r\subset \mathbb{T}^n$ un sous-tore dont l'action restreinte
associ\'ee sur $\mathcal{M}$ laisse invariantes les foncions $H_{\varepsilon}$, pour
tout $\varepsilon$. Notons par $i:\mathfrak{t}^r\hookrightarrow \mathfrak{t}^n$
l'inclusion de l'alg\`ebre de Lie de $\mathbb{T}^r$ dans celle de $\mathbb{T}^n$ et
par $i^*:(\mathfrak{t}^n) ^\ast \rightarrow (\mathfrak{t}^r) ^\ast $ son dual. Si
$m$ est $i^*\mu$-non-d\'eg\'en\'er\'e dans le sens de la
D\'efinition~\ref{nondegeneracy definition} et sa vitesse appartient \`a l'alg\`ebre
de Lie de $\mathbb{T}^r$ alors pour chaque valeur du param\`etre $\varepsilon$
suffisamment proche de $0$, les \'equilibres relatifs de $X_{H_{\varepsilon}}$ sont
en correspondance bijective avec les points critiques d'une fonction lisse
$\overline{h} _\varepsilon : \mathbb{T}^{n-r}\rightarrow \mathbb{R}$. Par
cons\'equent, il existe au moins $n-r+1$ \'equilibres relatifs de
$X_{H_{\varepsilon}}$ avec moment $i^*\mu$ et vitesse proche de $\xi$. De plus, si
les points critiques de $\overline{h} _\varepsilon$ sont tous non-d\'eg\'en\'er\'es
le nombre d'\'equilibres relatifs bifurqu\'es est au moins $2^{n-r}$.
\end{theoreme}

Pour \'etablir ce th\'eor\`eme, nous utilisons la technique de r\'eduction
symplectique par rapport au tore ${\mathbb T} ^r $ (espace de Marsden-Weinstein
${\mathcal M}_\alpha$) dans les coordonn\'ees de Marle-Guillemin-Sternberg
correspondant \`a l'action de $ \mathbb{T}^n $. Nous y caract\'erisons les
\'equilibres relatifs bifurquant de $m$ apr\`es brisure forc\'ee de sym\'etrie comme
points critiques de l'hamiltonien r\'eduit $h _\varepsilon$ sur  ${\mathcal
M}_\alpha$. Moyennant une condition de non-d\'eg\'en\'erescence sur l'\'equilibre
relatif initial $m$, nous obtenons de $h _\varepsilon$ une fonction $\overline{h}
_\varepsilon$ sur le tore ${\mathbb T}^{n-r}$ dont les points critiques sont en
correspondance bijective avec les \textsc{er} recherch\'es. Il s'en suit par la
th\'eorie des points critiques d'une fonction \`a variable r\'eelle sur une
vari\'et\'e compacte les estimations sur leur nombre (cat\'egorie de
Ljusternik-Schnirelmann, th\'eorie de Morse).

\selectlanguage{english}

\section{Introduction}

\begin{wrapfigure}{r}{100pt}
 \psfrag{a}{a}
 \psfrag{b}{b}
 \psfrag{mg}{mg}
\includegraphics[height=66pt,width=80pt]{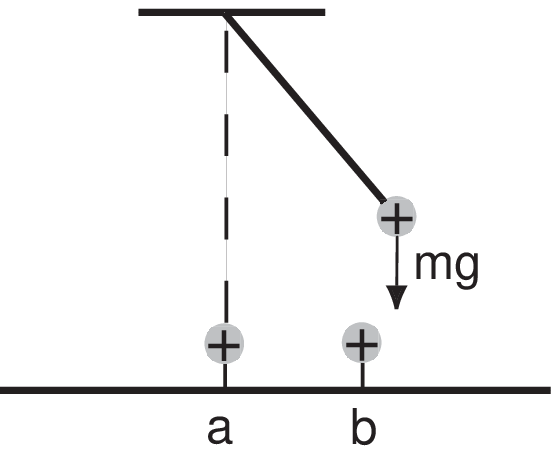}
\end{wrapfigure}
Forced symmetry breaking in dynamical systems is a phenomenon that takes place when
we add to a symmetric system a perturbation with less symmetry. In this note  we
study this phenomenon in the context of globally Hamiltonian dynamical systems, that
is, symmetric Hamiltonian systems to which a momentum map can be associated. The
particular problem of study is the ``survival'' or persistence of relative
equilibria of the fully symmetric system after the symmetry breaking perturbation is
added to it. Our motivation relies strongly upon the fact that this phenomenon is
naturally present in many systems. For instance, consider a spherical pendulum whose
bob of mass $m$ has been charged with a positive charge $q$ (see the figure).
Suppose now that right below the point of suspension of the pendulum we place a
charge identical to that of the pendulum (position (a) in the picture). If the
repulsive electrostatic force is strong enough, the stable downright equilibrium of
the spherical pendulum becomes unstable and a ring of equilibria appears. Suppose
now that the circular symmetry of the system is broken by slightly sliding the
charge to a side (position (b) in the picture). It can be seen that only two of the
equilibria in the ring survive. Our main goal in this paper is the formulation of a
general theorem capable of predicting such behaviour. The kind of systems we are
interested in can be mathematically described by considering a finite dimensional
symplectic manifold $({\mathcal M}, \omega )$ acted freely and canonically upon by
the $n$-torus ${{\mathbb T}}^n$. We assume that this action has a momentum map $
{\bf J}_{{{\mathbb T}}^n} : {\mathcal M} \longrightarrow ( \mathfrak{t}^n)
^\ast\simeq \mathbb{R} ^n $. Let $\mathbb{T} ^r\subset \mathbb{T} ^n $ be a
subtorus, $H_\varepsilon$  a family of Hamiltonian functions on ${\mathcal M}$
parametrized by $\varepsilon \in \mathbb{R} $, and assume that $H _0 $ is
$\mathbb{T}^n $-invariant whereas $H _\varepsilon $ is only $\mathbb{T}^r
$-invariant, for all $\varepsilon \in \mathbb{R}$. The problem that we discuss in
this note is under what conditions a given relative equilibrium $m \in  {\mathcal
M}$ of $H _0 $ with respect to its $\mathbb{T}^n  $-symmetry {\it persists} to
relative equilibria of the Hamiltonian vector fields associated to $H _\varepsilon
$, for $\varepsilon$ sufficiently small, with respect to their $\mathbb{T}^r
$-symmetry. We recall that a relative equilibrium (\textsc{re}) of a $\mathbb{T}^n
$-equivariant dynamical system $X $ on ${\mathcal M}  $ is a point $m$ for which
there exists an element $\xi $ in the Lie algebra $\mathfrak{t}^n  $ of $\mathbb{T}
^n$ (called the {\it velocity} of the \textsc{re}) such that $X (m)= \xi_{{\mathcal
M}}(m) $. The symbol $\xi_{{\mathcal M}}(m):= \left.\frac{d}{dt}\right|_{t=0} \exp t
\xi \cdot m $ denotes the infinitesimal generator of the $\mathbb{T}^n $-action
associated to the element $\xi$.

\section{Preliminaries}
Throughout we assume that $(\mathcal{M},\omega)$ is a finite-dimensional symplectic
manifold with a free Hamiltonian action of the torus ${\mathbb T}^n$ with momentum
map ${\mathbf J}_{{\mathbb T}^n}:\mathcal{M}\to( \mathfrak{t}^n) ^\ast$.  Such a
momentum map is necessarily \emph{invariant}: ${\mathbf J}(g\cdot m) = {\mathbf
J}(m)$ (with $g\in {\mathbb T}^n$ and $m\in\mathcal{M}$). We fix a torus subgroup
${\mathbb T}^r\subset {\mathbb T}^n$, and let
$i:\mathfrak{t}^r\hookrightarrow\mathfrak{t} ^n$ be the inclusion of Lie algebras.
The momentum map for the restricted action by ${\mathbb T}^r$ is ${\mathbf
J}_{{\mathbb T}^r} = i^*\circ{\mathbf J}_{{\mathbb T}^n}$, where
$i^*:(\mathfrak{t}^n)^*\rightarrow(\mathfrak{t}^r)^*$ is the dual map to $i$.

We also assume that $m\in \mathcal{M}$ is a relative equilibrium for the Hamiltonian
system $({\mathcal M}, \omega, H_0, {{\mathbb T}}^n, {\mathbf J}_{{{\mathbb
T}}^n})$, with momentum $\mu:={\mathbf J}_{{{\mathbb T}}^n}(m)\in (\mathfrak{t} ^n)
^\ast $ and velocity $\xi\in \mathfrak{t}^n $. We recall that this amounts to the
point $m$ being a critical point of the \emph{augmented Hamiltonian} $H _0-
\mathbf{J}_{\mathbb{T} ^n} ^\xi $, that is, $D(H_0-{\mathbf J}_{{{\mathbb
T}}^n}^{\xi})(m)=0$.

In order to formulate the main hypothesis of the theorem we need to recall the {\bf
Witt-Artin decomposition} of the tangent space $T_m{\mathcal M}$; define $V_m$ as a
complement in $\ker T_m{\mathbf J}_{{{\mathbb T}}^n}$ to the tangent space
${\mathfrak t}^n\cdot m$ at $m$ of the $\mathbb{T} ^n $-group orbit, that is, $ \ker
T_m{\mathbf J}_{{{\mathbb T}}^n}= V_m\oplus {\mathfrak t}^n\cdot m $. The
space
$V_m$ is called the {\bf symplectic normal space} at $m$. Notice that
${\mathfrak t}^n\cdot m
\subset (V _m)^{\omega(m)} $. Let $W$ be a Lagrangian complement to ${\mathfrak
t}^n\cdot m$ in $(V _m)^{\omega(m)}$. The decomposition
\begin{equation}
\label{witt artin decomposition} T_m{\mathcal M}=V_m\oplus  {\mathfrak
t}^n\cdot m\oplus W
\end{equation}
is called a Witt-Artin decomposion of the tangent space $T_m{\mathcal M}$. We will
refer to $W$ as the {\bf orbital complement} at $m$ of the Witt-Artin
decomposition~(\ref{witt artin decomposition}). To finish these preliminaries, we
give a definition which we will use in our result.

\begin{e-definition}
\label{nondegeneracy definition} With the notation as above, a {\bf nondegeneracy
space} ${\mathcal N}_\alpha $ at $m$ associated to the momentum $\alpha\in
(\mathfrak{t} ^r) ^\ast $ is defined as
\[ {\mathcal N}_{\alpha}=A_{\alpha}\oplus
V_m
\]
where $A_{\alpha}:=\left\{ w \in W\    \vert \   i^*(\mu +T_m{\mathbf J}_{{{\mathbb
T}}^n}(w))=\alpha \right\}$. Let $H \in C^{\infty}({\mathcal M})^{\mathbb{T}^n} $ be
a smooth $\mathbb{T} ^n $-invariant function on  ${\mathcal M} $ that exhibits a
critical point at $m$, that is, $DH (m)=0 $. We say that $m$ is an $\alpha$-{\bf
nondegenerate} critical point of $H$  when the symmetric bilinear form
\[
D^2H (m)|_{{\mathcal N}_\alpha\times {\mathcal N}_\alpha}
\]
is nondegenerate.
\end{e-definition}

\begin{lemma}
The $\alpha $-nondegeneracy of a critical point given in the previous definition
depends only on the value $\alpha\in (\mathfrak{t} ^r) ^\ast $ and not on the
specific Witt-Artin decomposition used to verify this condition.
\end{lemma}

\noindent\textbf{Proof.\ \ }It suffices to show that our nondegeneracy condition is
independent of the choice of $V _m$ and $W$ in the Witt-Artin decomposition.  Assume
that $H$ is $\alpha$-nondegenerate at $m$ for a fixed choice of $V _m$ and $W $. Let
$V _m' $ be another choice of symplectic normal space at $m$, $W' $ a complement to
${\mathfrak t}^n\cdot m $ in $(V_m')^{\omega(m)}$, and ${\mathcal N}_\alpha' $ the
associated nondegeneracy space. Let $v _1+ w _1, v _2 + w  _2 \in {\mathcal N}_
\alpha $ be arbitrary with $v _1, v _2 \in A _\alpha $ and $w _1, w _2 \in  V _m $.
The Witt-Artin decomposition of $T _m {\mathcal M} $ implies the existence of unique
elements $\xi, \eta \in \mathfrak{t} ^n $, $v _1', v _2' \in  A _\alpha' $, $w _1',
w _2' \in V _m ' $ such that
\[
v _1+ w _1= \xi _{{\mathcal M}} (m)+v _1'+ w _1' \quad\textrm{and}\quad
v _2+ w _2= \eta _{{\mathcal M}} (m)+v _2'+ w _2'.
\]
The $\mathbb{T}^n $-invariance of $H$ implies that
\[
D^2H (m)(v _1+ w _1,v _2+ w _2)=D^2H (m)(v _1'+ w _1',v _2'+ w
_2').
\]
Given that the map $v + w   \in {\mathcal N}_\alpha \longmapsto v' + w' \in {\mathcal
N}_\alpha' $ is an isomorphism, the result follows. \hfill $\Box$

\section{Theorem on forced symmetry breaking}

The goal of this section is to prove the following theorem:

\begin{theorem}
\label{main theorem on bifurcation symmetry breaking} Let  $({\mathcal M},
\omega,{\mathbb T}^n,{\mathbf J}_{{\mathbb T}^n}, {\mathbb T}^r, H_0,m,\xi,\mu )$ be as above. Let $H_\varepsilon$ be
a family of\/ ${\mathbb T}^r$-invariant Hamiltonian functions on ${\mathcal M}$ smoothly
parametrized by $\varepsilon \in \mathbb{R} $ with $H _0 $ $\mathbb{T}^n $-invariant.
Suppose that the relative equilibrium $m \in {\mathcal M} $ has velocity $\xi$. If
\begin{description}
\item  [(i)] $\xi \in i(\mathfrak{t}^r) $ and \item [(ii)] $m$ is a $i ^\ast \mu
$-nondegenerate critical point of $H _0- \mathbf{J} ^\xi _{\mathbb{T}^n} $, where
$\mu:= \mathbf{J}_{\mathbb{T} ^n} (m) $,
\end{description}
then for any value of the parameter $\varepsilon$ close enough to zero, the relative
equilibria of the Hamiltonian vector field $X _{H _\varepsilon} $ are in bijective
correspondence with the critical points of a smooth function
$\overline{h}_\varepsilon: \mathbb{T}^{n-r} \rightarrow  \mathbb{R} $. Consequently,
under these hypotheses there exist at least $n-r+1$  relative equilibria of
$X_{H_{\varepsilon}}$ with momentum $i^*\mu$ and velocity close to $\xi$.
Additionally, if the critical points of $\overline{h} _\varepsilon\in
C^{\infty}(\mathbb{T}^{n-r}) $ are all nondegenerate the number of bifurcated
relative equilibria is at least $2^{n-r}$.
\end{theorem}

\noindent\textbf{Proof.\ \ } The local character of the result that we want to
prove permits us to use the local model around the $\mathbb{T}^n $-orbit of $m$
given by the {\it Marle-Guillemin-Sternberg} (MGS) normal form
\cite{normal1,normal2}. We recall that  this result provides a $\mathbb{T}^n $-equivariant
symplectomorphism between a $\mathbb{T}^n $-invariant neighborhood of the orbit
$\mathbb{T}^n \cdot  m $ and the product $Y:=\mathbb{T}^n\times
(\mathfrak{t}^n) ^\ast
\times V _m $, considered as a $\mathbb{T}^n $-symplectic space with the
$\mathbb{T}^n $-action given by  $g \cdot (h, \eta, v):=(gh, \eta, v)  $, $g,h
\in \mathbb{T} ^n $, $ \eta \in (\mathfrak{t}^n) ^\ast$, $v \in V _m $
and with a symplectic form with respect to which the momentum map
associated to this
$\mathbb{T}^n $-action has the form ${\mathbf J}_{\mathbb{T}^n} (g, \eta , v)=
\mu+ \eta $. In this model, the point $m\in  {\mathcal M}$ is represented by
$(e,0,0)\in Y$ and the space $V _m $  is one of the symplectic normal spaces at
$m$ that we have previously defined. We will carry out the proof of our theorem
in these coordinates by looking for the critical points of the reduced
Hamiltonians  $h _{\varepsilon, \alpha} $ on the $ \mathbb{T} ^r $-Marsden-Weinstein
reduced space ${\mathcal M}_\alpha:= {\mathbf J} ^{-1}_{\mathbb{T} ^r}(\alpha)/
\mathbb{T}^r $ defined by $h _{\varepsilon, \alpha}\circ \pi_\alpha=H _\varepsilon \circ i
_\alpha $, where $\alpha:= i ^\ast  \mu $, $i _\alpha: {\mathbf J}
^{-1}_{\mathbb{T} ^r}(\alpha) \hookrightarrow {\mathcal M}  $ is the injection, and
$\pi _\alpha: {\mathbf J} ^{-1}_{\mathbb{T} ^r}(\alpha) \rightarrow {\mathcal
M}_\alpha $ is the projection. A straighforward computation  in MGS coordinates
shows that
\begin{equation}\label{eq4} {\mathbf J}^{-1}_{{{\mathbb T}}^r}(\alpha)={{\mathbb T}}^n\times
A_{\alpha}\times V_m
\end{equation}
where $A_{\alpha}$ is the vector subspace of $(\mathfrak{t}^n) ^\ast$ given by
$A_{\alpha}:=\left\{ \eta \in (\mathfrak{t}^n) ^\ast\ \vert \ i^*(\mu +\eta)=\alpha
\right\}=\ker i ^\ast  $ and that  ${\mathcal N}_\alpha:=A _\alpha \times V _m  $ is
a $\alpha$-nondegeneracy space at $m$. From   expression (\ref{eq4}) it is clear
that
\[
{\mathcal M}_{\alpha}={{\mathbb T}}^n\times A_{\alpha}\times
V_m/{{\mathbb T}}^r\simeq{{\mathbb T}}^{n-r}\times {\mathcal N}_{\alpha}.
\]
The problem of finding the relative equilibria in the statement of the theorem
is now equivalent to the search of the critical points of the real-valued
functions $h _{\varepsilon, \alpha} $ defined on the Marsden-Weinstein reduced space
${\mathcal M}_{\alpha}={{\mathbb T}}^{n-r}\times {\mathcal N}_{\alpha}$.

The hypothesis on the $\alpha$-nondegeneracy of $m$ as a critical point of $H _0-
\mathbf{J} ^\xi_{\mathbb{T} ^n} $ implies that the quadratic form $D^2{h_{0,
\alpha}}(e,0,0)|_{{\mathcal N} _\alpha \times  {\mathcal N}_\alpha}$ is
nondegenerate. In order to lighten the notation we will omit the symbol $\alpha $ in
the function $h _{\varepsilon, \alpha} $ in all that follows. With this notation, we
need to find the triples $(k, \tilde{v}, \varepsilon)\in {{\mathbb
T}}^{n-r}\times{\mathcal N}_{\alpha}\times{\mathbb R}$ such that
\begin{equation}\label{eq7}
  D h_{\varepsilon}(k,\tilde{v})=0
\end{equation}
We proceed by using the Implicit Function Theorem to eliminate the parameter
$\tilde{v}\in {\mathcal N}_{\alpha}$ from the equation (\ref{eq7}) by writing it in
terms of the $ {{\mathbb T}}^{n-r}$ and ${\mathbb R}$ variables. Indeed, consider the following
map
\[
\begin{array}{ccc} {\mathcal F} : {{\mathbb T}}^{n-r}\times{\mathcal N}_{\alpha}\times{\mathbb R}
&\longrightarrow & \left({\mathcal N}_{\alpha}\right)^*\simeq {\mathcal N}_{\alpha}\\
(k, \tilde{v},\varepsilon) &\longmapsto &D_{{\mathcal
N}_{\alpha}}h_{\varepsilon}(k,\tilde{v})\\
\end{array}
\]

Since $m\equiv (e,0,0)$ is a $\mathbb{T} ^n $-relative equilibrium for $H_0$ we
have  ${\mathcal F} (g,0,0)=0$, for all $g\in {{\mathbb T}}^{n-r}$. Moreover, since the
partial derivative $D_{{\mathcal N}_{\alpha}}{\mathcal F}(g,0,0) : {\mathcal
N}_{\alpha}\longrightarrow \left({\mathcal N}_{\alpha}\right)^*\simeq {\mathcal N}_{\alpha}$ of
${\mathcal F}$ with respect to the ${\mathcal N}_{\alpha}$-factor, evaluated at
$(g,0,0)$ is given by  $D_{{\mathcal N}_{\alpha}}{\mathcal F}(g,0,0)=D^2_{{\mathcal
N}_{\alpha}}h_0(g,0)$ then the hypothesis on the $\alpha$-nondegeneracy of $m$
as a critical point of $H _0- \mathbf{J} ^\xi_{\mathbb{T} ^n} $ implies that
$D_{{\mathcal N}_{\alpha}}{\mathcal F}(g,0,0) : {\mathcal N}_{\alpha}\longrightarrow {\mathcal
N}_{\alpha}$  is injective. Consequently, $D_{{\mathcal N}_{\alpha}}{\mathcal
F}(g,0,0)$ is an isomorphism and we can then define via  the Implicit Function
Theorem  a smooth map $\tilde{v}_g: {\mathcal U}_{g}\times {\mathcal W}_g \longrightarrow
\left({\mathcal N}_{\alpha}\right)_g\subset {\mathcal N}_{\alpha}$ defined in an open
neighborhood of $(g,0)\in {{\mathbb T}}^{n-r}\times{\mathbb R}$ such that, for any $(k,\varepsilon)\in
{{\mathbb T}}^{n-r}\times{\mathbb R}$ in that neighborhood, we have that:
\begin{equation}\label{eqv}
  D_{{\mathcal N}_{\alpha}}h_{\varepsilon}(k,\tilde{v}_g(k,\varepsilon))=D_{{\mathcal
 N}_{\alpha}}h_{\varepsilon}(k,{\tilde{v}}_{g,\varepsilon}(k))=0
\end{equation}
Given that this argument can be repeated for any $g \in \mathbb{T}^{n-r}$ we can
invoke the compactness of $\mathbb{T}^{n-r} $ to build a finite family of functions
${\tilde{v}}_{g_i}:{\mathcal U}_{g_i}\times {\mathcal W}_{g_i} \longrightarrow
\left({\mathcal N}_{\alpha}\right)_{g_i}$, $i \in \{1, \ldots, \ell\} $, satisfying
(\ref{eqv}) and such that   $\bigcup _{i=1}^\ell{\mathcal U}_{g_i}=
\mathbb{T}^{n-r}$. Let us define

\[
\begin{array}{ccc} {\tilde{v}} : {{\mathbb T}}^{n-r}\times
\bigcap_{i=1}^\ell{\mathcal W}_{g_i}  &\longrightarrow & \bigcup_{i=1}^\ell
\left({\mathcal N}_{\alpha}\right)_{g_i}\\ (g,\varepsilon) &\longmapsto
&{\tilde{v}}_{g _i}(g,\varepsilon) \
\       \mbox{if }    g\in {\mathcal U}_{g_i}.\\
\end{array}
\]
This map is well defined by  the uniqueness of the maps ${\tilde{v}}_{g _i}$
obtained from the Implicit Function Theorem. Taking into account this new map,
our bifurcation equation  (\ref{eq7}) is now equivalent to:
\begin{equation}\label{eq9}
Dh_{\varepsilon}(k,{\tilde{v}}(k, \varepsilon))=0.
\end{equation}
The solutions of this equation coincide with the critical points of the
function ${\overline h}_{\varepsilon}(k):=h_{\varepsilon}(k,{\tilde{v}}(k, \varepsilon))$ defined,
for each value of the parameter $\varepsilon$,  on the compact manifold ${{\mathbb T}}^{n-r}$.
Indeed, using~(\ref{eqv}), we have
\begin{eqnarray}
D{\overline
h}_{\varepsilon}(t)&=&D_{{{\mathbb T}}^{n-r}}
h_{\varepsilon}(t, \tilde{v}(t,\varepsilon)) +D_{{\mathcal N}_{\alpha}}h_{\varepsilon}(t,
\tilde{v}(t,\varepsilon))\cdot D_{{{\mathbb T}}^{n-r}}\tilde{v}(t,\varepsilon))=D_{{{\mathbb T}}^{n-r}}
h_{\varepsilon}(t, \tilde{v}(t,\varepsilon)) \nonumber \\
    &=&D_{{{\mathbb T}}^{n-r}}
h_{\varepsilon}(t, \tilde{v}(t,\varepsilon))+D_{{\mathcal N}_{\alpha}}h_{\varepsilon}(t,
\tilde{v}(t,\varepsilon))=Dh_{\varepsilon}(t,{\tilde{v}}(t, \varepsilon)).
\nonumber
\end{eqnarray}
Consequently, the pair $(t,{\tilde{v}}(t, \varepsilon))$ is a solution
of~(\ref{eq9}) if and only if $t\in {{\mathbb T}}^{n-r}$ is a critical point of
${\overline h}_{\varepsilon} $. A lower bound for the number of these critical
points is provided by the Ljusternik-Schnirelmann category   ${\rm Cat}({{\mathbb
T}}^{n-r})=n-r+1$ of the torus ${{\mathbb T}}^{n-r} $ (see for instance
\cite{Ljusternik 1966}), which proves the statement of the theorem. Additionally, if
we know in advance that the critical points of ${\overline h}_{\varepsilon} \in
C^{\infty}(\mathbb{T}^{n-r})$ are all nondegenerate, the Morse inequalities
guarantee that this function has at least
$b^0(\mathbb{T}^{n-r})+b^1(\mathbb{T}^{n-r})+ \cdots+b^{n-r}(\mathbb{T}^{n-r})$
critical points, where $b^i (\mathbb{T}^{n-r}) $, $i \in \{0, \ldots, n-r \} $, is
the $i$-th Betti number of the torus $\mathbb{T}^{n-r} $. Since $b^i
(\mathbb{T}^{n-r}) = \pmatrix{n-r\cr i} $, $i \in \{0, \ldots, n-r \} $, we have
\[
b^0(\mathbb{T}^{n-r})+b^1(\mathbb{T}^{n-r})+ \cdots+b^{n-r}(\mathbb{T}^{n-r})=\sum
_{i=0}^{n-r} \pmatrix{n-r\cr i} =2^{n-r}
\]
and hence the second estimate in the statement follows. \hfill $\Box$

\section{Symmetry breaking using Poisson reduction}

In the previous theorem we confined our search for bifurcated relative equilibria to
the  momentum level set $\mathbf{J}_{\mathbb{T}^r}^{-1}(i ^\ast\mu)$. This fact
appears in the proof of that result when we use the symplectic reduced space
${\mathcal M}_{i ^\ast \mu} $. If instead of using  ${\mathcal M}_{i ^\ast \mu} $ we
consider the Poisson reduced space $\widetilde{{\mathcal M}}:={\mathcal M}/
\mathbb{T}^r $ we can obtain another bifurcation result where the predicted relative
equilibria could, in principle, have a momentum different from that of the given
\textsc{re}. This is obtained at the expense of imposing a more demanding
nondegeneracy condition.

\begin{theorem} \label{poisson}
Let  $({\mathcal M}, \omega,{\mathbb T}^n,{\mathbf J}_{{\mathbb T}^n}, {\mathbb
T}^r, H_0,m,\xi,\mu )$ be as in Section 2. Let  $H_\varepsilon$  be a family of
Hamiltonian functions on ${\mathcal M}$ parametrized by $\varepsilon \in \mathbb{R}
$ and assume that $H _0 $ is $\mathbb{T}^n $-invariant whereas $H _\varepsilon $ is
only $\mathbb{T}^r $-invariant, for all $\varepsilon \in \mathbb{R}$. Suppose that
the point $m \in {\mathcal M} $ is a $ \mathbb{T}^n $-relative equilibrium of the
Hamiltonian vector field $X_{H _0}$ with velocity $\xi\in i(\mathfrak{t}^r)$.
Suppose moreover that
\[
D^2\left.(H _0- \mathbf{J} ^\xi _{\mathbb{T}^n})(m) \right\vert_{{\mathcal N}\times
{\mathcal N}}
\]
is a nondegenerate quadratic form, where ${\mathcal N}:=W\times V_m $, for the
symplectic normal space $V _m $ and orbital complement $W $ corresponding to some
Witt-Artin decomposition of $T _m {\mathcal M} $. Then for any value of the
parameter $\varepsilon$ close enough to zero, the relative equilibria of the
Hamiltonian vector field $X _{H _\varepsilon} $ are in bijective correspondence with
the critical points of a smooth function $[h_{\varepsilon}] : {{\mathbb
T}}^{n-r}\longrightarrow {\mathbb R}$. Consequently, under these hypotheses there
exist at least $n-r+1$  $\mathbb{T}^r $-relative equilibria  near $m$ with momentum
close to $i^*\mu$ and velocity close to $\xi$.  Additionally, if the critical points
of $[h_{\varepsilon}]\in C^{\infty}(\mathbb{T}^{n-r}) $ are all nondegenerate the
number of these bifurcated relative equilibria is at least $2^{n-r}$.
\end{theorem}

\noindent\textbf{Proof.\ \ } This mimics the proof of Theorem~\ref{main theorem on
bifurcation symmetry breaking} where the reduced space ${\mathcal M} _\alpha $ has
been replaced by ${\widetilde{\mathcal M}} $. Note that in the MGS normal form
coordinates we can  write, locally,
\[
{\widetilde{\mathcal M}}=\left({{\mathbb T}}^n\times (\mathfrak{t} ^n) ^\ast \times
V_m\right)/{{\mathbb T}}^r\simeq {{\mathbb T}}^{n-r}\times (\mathfrak{t} ^n)
^\ast\times V_m.
\]

\vskip -6mm

\hfill$\Box$

\section*{Acknowledgements}
This
research was partially supported by the European Commission through
funding for the Research Training Network
\emph{Mechanics and Symmetry in Europe} (MASIE).

\end{document}